\documentclass[12pt,leqno,]{article}
\usepackage{latexsym}
\usepackage[all]{xy}
\usepackage{theorem, latexsym, amsfonts}
\usepackage{amssymb, amsmath}
\usepackage{color}
\newcommand{\beq}{\begin{equation}}
\newcommand{\eeq}{\end{equation}}

\begin{document}

\title{\textbf{A Model Categorical Approach to Group Completion of $\mathbf{E_n}$-Algebras }}
\author{M. Stelzer }

\date{\today}
\maketitle

\noindent\begin{it}Abstract: A group completion functor $Q$ is constructed in the category of algebras in simplicial sets over a cofibrant $E_n$-operad $\mathcal{M}$. It is shown that $Q$ defines a Bousfield-Friedlander simplicial model category on $\mathcal{M}$-algebras.\end{it}

\noindent\textit{Keywords: 2010 MSC: 55P35, 55P48, 55P60.}

\setcounter{section}{-1}
\section{Introduction} 

A space $X$ is of the homotopy type of an n-fold loop space if and only if it carries an action of an $E_n$-operad which induces a group-like H-structure on $X$.
 Since algebras over a suitable
 $E_n$-operad  form a 
model category \cite{SV} it is natural to consider the localization of these model structures
 with respect to some variant of group completion functor $Q$. 
 The aim of this paper is to construct such a Bousfield-Friedlander $Q$-structure for algebras over a cofibrant
 $E_n$-operad $\mathcal{M}$ in simplicial  sets. The fibrant objects turn out to be essentially the group complete objects whose underlying simplicial sets are fibrant. Hence, in the light of the delooping results in \cite{BV}, \cite{M}, the  $Q$-local homotopy category may be viewed as the homotopy category of n-fold loop spaces.
The results of this paper play a role  in our current  joint work \cite{FSV} with Zig Fiedorowicz
and Rainer Vogt on n-fold monoidal categories \cite{B}.
		
The paper is organised as follows. In section 1, we construct a group completion functor $\bar{Q}$ together with a coaugmentation  $\bar{q}:1\to \bar{Q}$ in topological algebras over a cofibrant $E_n$-operad. For technical reasons, we use  a mixture of classical and model categorical homotopy theory. Since an $E_n$-space is in general not a monoid we have to change operads. For this step,  we rely on a recent theorem of Morton Brun, Zig Fiedorowicz
and Rainer Vogt on tensor products of operads \cite{BFV}. This result is used to generalize and adapt an argument of May, from $n=\infty$  to all $n$.
  In the next section,  we recall Bousfield's 
improved axioms for the existence of a localized model structure defined by a coaugmented functor $Q$.
It is  shown that a $Q$-local model structure exists for algebras over a cofibrant $E_n$-operad $\mathcal{M}$  in simplicial sets. The functor $Q$ is induced from the topological version $\bar{Q}$ defined in  section 1.  In a nutshell, the strategy for the proof of the main result Theorem 2.7. may be described as follows. In order to verify Bousfield's axioms we have to use properties of the classifying space functor $B$ which figures in the classical group completion $\Omega B M$ of a topological monoid $M$. On the other hand, to construct a coaugmented functor in $E_n$-algebras we have to rely on May's machine and on the cofibrancy of the operad  $\mathcal{M}$. The point is, that we need a natural transformation $\bar{q}$ which commutes with the operad action on the nose and not only up to coherent homotopy. So we have to merge the best of two worlds in our construction. 

In \cite{BCV} the authors establish a related result. They consider the model category of algebras over the theory which 
encodes all natural maps between products of n-fold loop spaces and prove a recognition theorem for n-fold loop spaces 
with respect to this theory. Also, in unpublished work, Bousfield generalizes Segal's approach \cite{S} to infinite loop space theory to n-fold loop
 spaces \cite{BO2}. 
However,  for the applications which we have in mind, the close connection to the operad of little n-cubes is essential. 
We will use freely the language of model categories. Besides the original source \cite{Q} there are now
some more recent  books on this subject  \cite{GJ}, \cite{HI}, \cite{H}.
For general background information
on operads the reader may consult \cite{MSS}. 

In this paper, we will work in the ground categories $Top$, $Top_{\ast}$, $SS$, $SS_{\ast}$ of k-spaces, based k-spaces  \cite{V1}, simplicial sets and based simplicial sets. For an operad $M$ in $Top$ or $SS$ we write $Top_{\ast}^M$ and $SS_{\ast}^M$  for the categories of $M$-algebras.

\noindent{\bf Acknowledgements:} I would like to thank Oliver Roendigs for help in the proof of Theorem 1.6. and Rainer Vogt
for discussions and for spotting an error in an earlier version of the paper.

\section{Group completions of $E_{n}$-algebras}

The group completion of a topological monoid $M$ is the loop space of the classifying space $\Omega BM$. We use Milgram's version of the classifying space construction in this paper  \cite{MI}. There is a natural map 
\[M\stackrel{\kappa}{\to}\Omega BM\] which is well known to be an $A_{\infty}$-map \cite{BV}. This means $\kappa$ commutes with the action of an $A_{\infty}$-operad up to coherent homotopy.
If $X$ is an $A_{\infty}$-space there is a functor \cite{BV} which replaces $X$ by an equivalent monoid $MX$ 
and the group completion of $X$ can be defined as $\Omega BMX$.
 For H-spaces whose multiplication satisfies a weak form of homotopy commutativity there is a   homological version of group completion which we recall in the form given in \cite{M2}.

\noindent{\bf Definition 1.1. }\begin{it} An H-space $X$ is called admissible if $X$ is homotopy associative and if left translation by any given element of 
$X$ is homotopic to right translation with the same element. An H-map between admissible H-spaces 
\[g:X\to Y \]
is called a homological group completion of $X$ if $Y$ is group-like and the unique morphism of $k$-algebras
\[\bar{g}_{\ast}:H_{\ast}(X;k)\left[\pi_0^{-1}X\right]\to H_{\ast}(Y;k)      \]
which extends $g_{\ast}$ is an isomorphism for all commutative coefficient rings
$k$.\end{it} 

The group completion theorem \cite{SMS}, \cite{QFM}(see also \cite{M2}) asserts that $\kappa$  
is a homological group completion for  a well pointed topological monoid  $M$ in case the multiplications of $M$ and $\Omega BM$ are admissible.  

 Our goal in this section is to construct a coaugmented  functor $\bar{Q}$ in topological $E_n$-algebras which is closely related to a classical group completion. It will turn out that this functor is a homological group completion for $n>1$. 
Many of our arguments are adaptions of the ones given by May for the case $n=\infty $ in \cite{M2}, complemented by
model categorical considerations, which are needed to ensure that the natural maps\[\bar{q}:X\to \bar{Q}X\] are $E_n$-homomorphisms.

Next we recall some notions and results from the literature.

\noindent{\bf Definition 1.2. }\begin{it} A topological operad $M$ is called:\\
$(i)$   well-pointed if the inclusion 
\[\{ id  \} \to M_1\]
is a closed cofibration,\\
$(ii)$ reduced if $M_0 =\ast$,\\
$(iii)$ $\Sigma$-free  if $M_{k}$ is $\Sigma_{k} $-free and $M_{k}\to M_{k}/\Sigma_{k} $ is a
 numerable principal $\Sigma_{k}$-bundle \cite{BV},\\
$(iv)$ $E_n$-operad if there exists a chain of maps of operads 
\[M=B^0 \stackrel{f^0}{\to}B^1 \stackrel{f^1}{\leftarrow}\cdots \to B^r \stackrel{f^r}{\leftarrow}\mathcal{C}_n \] such that
\[ f_i^s :B_i^s \to B_{i}^{s+1} \]  are $\Sigma_i$-equivariant homotopy equivalences, and where $\mathcal{C}_n$ is the operad of little $n$-cubes \cite{BV}.\\
A topological or simplicial operad $M$ is called:\\
$(v$) $\Sigma$-cofibrant if the underlying collection is cofibrant in the model category of collections \cite{BM1}.\\
$(vi)$ cofibrant if it is cofibrant in the model category of operads \cite{BM1}.
 \end{it}

Note that $\mathcal{C}_n$ is $\Sigma$-free but not known to be $\Sigma$-cofibrant.
In general, notions of cofibrancy related to classical homotopy theory, whose model categorical incarnation in $Top$
 is the Str\o m structure \cite{ST}, are weaker than the ones related to the Quillen model structure \cite{Q}.
 There are  natural ways to replace a given operad
 by a well-pointed  \cite{V3} or a cofibrant one.
For example, the Boardman-Vogt $W$-construction \cite{BV} serves as a cofibrant replacement functor.
This is shown in \cite{V3} in the setting of the cofibration category of topological operads
 with underlying Str\o m structure. The model categorical case is treated in \cite{BM2}.
There is also a reduced version of $W$ \cite{BV}, \cite{BM2} which is cofibrant in the model category of reduced operads.
 In the following, we will restrict attention to reduced operads.
Algebras over reduced operads have only one 0-ary operation.
 In the based context this operation is assumed to coincide with the basepoint.
In this situation, the relevant monad associated with a reduced operad
  is defined by certain identifications related to the base points \cite{M}.

 Let $C_n$ be the monad defined by $\mathcal{C}_n$ on pointed spaces.
 There is  a morphism of monads \cite[5.2.]{M} \[\alpha_n :C_n \to \Omega^{n}\Sigma^{n}.\]
such that the induced map of $\mathcal{C}_n$-algebras \[C_{n}X\to \Omega^{n}\Sigma^{n}X\] is a homological group completion for $n>1$ as has been shown by
 Cohen and Segal  \cite{C},\cite{S2}.

 In the proof of the proposition below will make use of May's
 two-sided bar construction \cite{M} and of the results in \cite[Appendix]{M}, \cite[Appendix]{M2}.

\noindent{\bf Proposition 1.3. }\begin{it}Let $M$ be a cofibrant, reduced topological $E_n$-operad. Then there is a functor 
\[G:Top_{\ast}^ {M}\to Top_{\ast}\] such that $GX$ is a topological monoid,
and a natural transformation $g$, with values in H-maps, from the forgetfull functor 
\[U:Top_{\ast}^ {M}\to Top_{\ast}\]
to $G$ such that 
\[g:X\to GX\]
is a homological group completion for well pointed $X$ in case $n>1$.\end{it}

\noindent{\bf Proof:} In a first step we change operads since we need a monoid structure to go along the algebra structure.
Let $\mathcal{A}$ be the operad whose algebras are topological monoids and $\mathcal{A}\otimes \mathcal C_{n-1}$
 the  tensor product of operads \cite{BV}. There is an operad map
 \[\gamma :M\to \mathcal{A}\otimes \mathcal C_{n-1}  \] whose underlying maps are homotopy equivalences.
This is the case since $M$ is cofibrant,
$\mathcal{A}\otimes \mathcal C_{n-1}$ is an $E_n$-operad  \cite[Theorem C]{BFV} and because   $\mathcal C_{n}$ is of the homotopy type of a CW-complex \cite{B}.
The operad  $\mathcal{A}\otimes \mathcal C_{n-1}$ is well-pointed and $\Sigma$-free \cite{BFV}.
 Define $G$  by
 \[GX=\Omega_M BB(\mathcal{A}\otimes \mathcal C_{n-1},M,X)\]
where $B(\mathcal{A}\otimes \mathcal C_{n-1},M,X)$ is the two sided bar construction, $\Omega_M$ the Moore loop functor, and define $g$ as the composition
\[ X\stackrel{\tau}{\to}B(M,M,X)\stackrel{B(\gamma ,1,1)}{\longrightarrow}B(\mathcal{A}\otimes \mathcal C_{n-1},M,X)\]
\[B(\mathcal{A}\otimes \mathcal C_{n-1},M,X)\stackrel{\bar{\kappa}}{\to}\Omega_M B B(\mathcal{A}\otimes \mathcal C_{n-1},M,X).\]
Here  $\bar{\kappa}$ is the natural inclusion and $\tau $ is the right inverse of the augmentation 
 \[X\stackrel{\epsilon}{\leftarrow}B(M,M,X)\] which is a $M$-homomorphism and strong deformation retraction.
 So $\tau $ is a homotopy equivalence and a $M$-map in the sense of \cite{BV}. In particular it is an H-map.
 Note that
  $B(\mathcal{A}\otimes \mathcal C_{n-1},M,X)$ is an  $\mathcal{A}\otimes \mathcal C_{n-1}$-space
 and hence it is a monoid in $\mathcal C_{n-1}$-spaces.

 The $M$-morphism  $B(\gamma ,1,1)$ is a homotopy equivalence for 
 well pointed $X$.   In order to see this we consider
 the map induced by $\gamma$ between the simplicial spaces whose realizations are the bar constructions in question
\[ B_{\ast}(M,M,X)\stackrel{B(\gamma ,1,1)_{\ast}}{\longrightarrow}B_{\ast}(\mathcal{A}\otimes \mathcal C_{n-1},M,X). \] These simplicial spaces are proper \cite{M}. To see this we apply  \cite[A.10.]{M}.
The assumptions made there hold since by A.7. loc.cit.  the monad defined by
 $\mathcal{A}\otimes \mathcal C_{n-1}$ is an admissible $M$-functor. We apply  in \cite[A.4.]{M2} which states that a map between proper simplicial spaces 
 which is a homotopy
equivalence in any given simplicial degree induces a homotopy equivalence after realization. 
 Since $M$ and $\mathcal{A}\otimes \mathcal C_{n-1}$ are both $\Sigma$-free the maps $B(\gamma ,1,1)_{m}$  are indeed homotopy equivalences \cite[A.2.]{M2}, \cite[A.3.4.]{BV}.\\
So the composition of the first two maps whose composition make up $g$ is a homotopy equivalence and an H-map between 
admissible H-spaces. The last assertion follows from the group completion theorem applied
to the third map since, using well known facts, one can replace $\Omega $ by $\Omega_M $.$\Box$

\noindent{\bf Proposition 1.4.} \begin{it} Let $M$ be a cofibrant, reduced topological $E_n$-operad and $X$ a well pointed $M$-algebra.
Then 
\begin{enumerate}
\item[(i)] the composition
\[X\stackrel{\tau}{\to}B(M,M,X)\stackrel{B(\alpha_n\pi_n ,1,1)}{\longrightarrow}B(\Omega^n \Sigma^n  ,M,X)\]
is a homological group completion for $n>1$.
\item[(ii)] the spaces $GX$ and $B(\Omega^n \Sigma^n  ,M,X)$ are naturally  weakly equivalent by   H-maps
for all $n>0$.
\end{enumerate} 
\end{it} 

\noindent{\bf Sketch of Proof:} The assertions in (i) and (ii) in the case $n>1$ follow by the same argument as in Theorem 2.3.(ii) in \cite{M2}, once one replaces Theorem 2.2 in loc.cit. by the result of Cohen and Segal mentioned above.
To settle the case   $n=1$, we appeal to Thomason's result on the uniqueness  of delooping machines \cite{T}, which  implies that the May delooping $B( \Sigma  ,M,X)$ and $BB(\mathcal{A}\otimes \mathcal C_{n-1},M,X)$ are homotopy equivalent for well pointed $X$. Since there is a  weak equivalence \[\rho:B(\Omega \Sigma  ,M,X)\to \Omega B( \Sigma  ,M,X)\] \cite[13.1(iii)]{M}  the assertion for $n=1$ is proved.
\hfill$\Box$

\noindent{\bf Definition 1.5.  }\begin{it} A model category $\mathcal{D}$ is called right proper if every pullback of a weak equivalence along a fibration is a weak equivalence, left proper if every pushout of a weak equivalence along a cofibration is a weak equivalence, and proper if it is right and left proper.\end{it}

As shown in \cite{SV}, \cite{SP}, \cite{BM1} the model structures on $Top_{\ast}^M$ and $SS_{\ast}^M$, for  a $\Sigma$-cofibrant operad $M$,  are transferred \cite{Cr} along the free $M$-algebra functor from
$Top_{\ast}$ and $SS_{\ast}$ and are cofibrantly generated. Moreover, it was shown by Spitzweck \cite[Theorem 4 in section 4]{SP}  that,
for a cofibrant operad $M$, the transferred structure  is right proper and  the pushout in $Top_{\ast}^M$ of a weak equivalence along a cofibration is a weak equivalence, provided the source is cofibrant in the underlying model category of $Top_{\ast}$. In $SS_{\ast}^M$ this holds unconditionally. Hence, the model category of algebras over a cofibrant operad in simplicial sets is proper.

Now we turn to the construction of the functor $\bar{Q}$:

\noindent{\bf Theorem 1.6.  }\begin{it} Let $\mathcal{M}$ be a reduced cofibrant simplicial operad such that the topological realization $M=|\mathcal{M}|$
 is an $E_n$-operad. Then there is  a  functor \[\bar{Q}:Top_{\ast}^ {M}\to Top_{\ast}^ {M}\]
and a natural transformation \[\bar{q}:1\to \bar{Q} \] in $Top_{\ast}^ {M}$ which is a homological group completion for algebras which are cofibrant  in  $Top_{\ast}$ if $n>1$.   
If $X\in Top_{\ast}^ {M}$ is  cofibrant in $Top_{\ast}$, then $\bar{Q}X$ is so as well and is naturally weakly equivalent to $GX$ by H-maps for all $n$. 
  \end{it}  

\noindent{\bf Proof:} Let $X$ be an $M$-algebra which is cofibrant as a pointed space.
Note that $X$  and $M$ are both well pointed  the later because 
geometric realization respects cofibrations.  
 The operad $M$ is cofibrant since $\mathcal{M}$ is and it is $\Sigma$-cofibrant by  \cite[4.3.]{BM1} and consequently also $\Sigma$-free \cite[Appendix 3.]{BV}.
Hence, $B(\alpha_n\pi_n ,1,1)$ is a homological group completion if $n>1$ by 1.4. Put \[HX=B(M ,M,X)\] and \[KX=B(\Omega^n \Sigma^n  ,M,X)\] for short.
These functors come with natural transformations of $M$-algebras  \[\epsilon :H \to 1\] and  
\[\eta = B(\alpha_n\pi_n ,1,1):H\to K\] where the underlying map of $\epsilon $ is a homotopy equivalence. Apply the natural CW-approximation $T=|S|$ to the simplicial spaces
whose realization is the diagram
\[X\stackrel{\epsilon}{\leftarrow} HX\stackrel{\eta}{\to} KX \] and realize the diagram to \[TX=X_T \stackrel{\epsilon_T}{\leftarrow} H_T X\stackrel{\eta_T}{\to}K_T X .\]
By \cite[4.9.]{FV} $TZ$ is a $M$-algebra with cellular action and the natural map $\phi :TZ\to Z$ is a weak equivalence of
 $M$-algebras for every $M$-algebra $Z$. It follows from \cite[11.4.]{M} that 
\[X\stackrel{\phi}{\leftarrow} TX\stackrel{\epsilon_T}{\leftarrow} H_T X\stackrel{\eta_T}{\to} K_T X \] is a diagram of $M$-algebras with all underlying spaces cofibrant. 

Let \[H_T X\stackrel{i}{\to}U\stackrel{j}{\to}K_T X\] be a natural factorization of $\eta_{T}$ into  cofibration $i$  and trivial
 fibration $j$ in the model category of $M$-algebras.

 Define $\bar{Q}X$ to be the pushout in $M$-algebras of  $i ,\phi \epsilon_{T}$ 
$$
\xymatrix{
H_{T}X \ar[r]^{\phi\epsilon_T} \ar[d]_i & 
X \ar[d]^{\bar{q}}
\\
U \ar[r]_h & \bar{Q}X
}
$$ 
   
 By \cite[Theorem 4 in section 4]{SP}, $h$ is a  weak equivalence 
 hence $\bar{Q}X$ with the induced natural transformation of $M$-algebras $\bar{q}: 1\to \bar{Q}$ is a group completion of
 $X$  in the category of $M$-algebras. The space $\bar{Q}X$ is cofibrant because the induced map $\bar{q}:X\to \bar{Q}X$ is a 
cofibration
 of algebras since $i$ is one. But a cofibration of algebras is a cofibration of spaces by \cite[Theorem 4 in section 4 ]{SP}. Now $X$ was assumed 
to be a cofibrant space hence $\bar{Q}X$ is one as well. \hfill$\Box$

\noindent{\bf Remark 1.7.} The assumption that $M$ is the realization of a simplicial operad is not severe because there is a Quillen equivalence between simplicial and topological operads \cite{BM1}.

\section{The $Q$-structure on $\mathcal{M}$-algebras }

 Let $\mathcal{D}$ be a  proper model category
and $Q:\mathcal{D} \to \mathcal{D}$ a coaugmented functor. 

Following Bousfield and Friedlander, we say that a morphism $f:X\to Y$ in $\mathcal{D}$  is a $Q$-equivalence if $Qf$ is a weak equivalence, a $Q$-cofibration if $f$
 is a cofibration, and a $Q$-fibration if $f$ satisfies the right lifting property with respect to $Q$-trivial 
 cofibrations.

Consider the following axioms: 

(A1) for each weak equivalence $f:X\to Y$  in  $\mathcal{D}$  
 the map \[Qf:QX\to QY\] is a weak equivalence; 

(A2) for each object in  $X\in\mathcal{D}$  the maps \[q_{QX},Qq_{X}
:QX\to QQX\] are weak equivalences;

(A3) for each pullback square  
$$
\xymatrix{
V \ar[r]^{k} \ar[d]_g & 
X \ar[d]^f
\\
W \ar[r]_h & Y
}
$$ 
in $\mathcal{D}$, with $f$ a fibration of fibrant objects such that
  $q:X\to QX$, $q:Y\to QY$ and $Qh:QW\to QY$ are  weak equivalences, the map 
\[QV\stackrel{Qk}{\to} QX \] is a weak equivalence.

It is a theorem of Bousfield \cite[9.3.]{BO} that in case (A1)-(A3) hold, then the three classes of maps given above define a proper model  category on $\mathcal{D}$. We will apply this theorem, or more precisely its  proof, to the category $SS_{\ast}^{\mathcal{M}}$ where $\mathcal{M}$ a cofibrant operad whose realization is an $E_n$-operad.
 As for the stable model structure of $\Gamma$-spaces constructed in \cite{BF} axiom (A3) does not hold in full generality.
So we have to adapt the arguments in \cite{BF} to the situation at hand. There is a $Q$-structure even for non proper $\mathcal{D}$. This follows from \cite[9.5]{BO} and \cite{STa}. However, we had to rely on the weak form of left properness of $Top_{\ast}^M$ for the proof of 1.6.

\noindent{\bf Lemma 2.1.}\begin{it} Let $M=|\mathcal{M}|$ be as in 1.6. and $\bar{Q}, \bar{q}$ the coaugmented functor constructed
 in section 1.
Then  the pair $\bar{Q}, \bar{q}$ satisfies (A1) and (A2) if $X$ and $Y$ are cofibrant spaces.\end{it} 

\noindent{\bf Proof:} By 1.6., we may replace $\bar{Q}X$ by $GX=\Omega_M BB(\mathcal{A}\otimes \mathcal C_{n-1},M,X)$ in the argument. Then the assertion  follows from some well known properties of
$\Omega_M$ and $B$ (see \cite[Chapter VI.]{BV}). Property (A1) is satisfied since the functors  $\Omega_M$ and $B$ preserve weak equivalences.
 For a  connected  space $X$ of the homotopy type of a CW-complex the natural evaluation map
\[e:B\Omega_M X\to X\]
\[e(t_1 ,x_1 ,\ldots ,t_k ,x_k )=
\omega (\sum_{i=1}^k (1-t_{1}\ast t_2 \ast \ldots \ast t_i )a_i)\]
is a homotopy equivalence \cite[6.15]{BV}
where $x_i =(\omega_i ,a_i ) \in \Omega_M X $, $t_1 \ast t_2 =t_1 +t_2 -t_1 t_2$, and $(\omega ,\sum_{i=1}a_i )=x_1\cdot x_2\cdot\ldots \cdot x_k $ 
.
The functor $\Omega_M B$ is a monad up to homotopy with structure morphisms $\Omega_M e$ and $\bar{\kappa}$. In particular $(\Omega_M e)( \bar{\kappa}_{\Omega_M BX})$ and $(\Omega_M e)(\Omega_M \bar{\kappa}_X )$ are homotopic to the identity. Property (A2) follows.  \hfill$\Box$

Recall that the adjoint pair $|-|, S$ geometric realization and the singular functor induces an adjoint pair \cite{FV} which will be denoted by
the same symbols $|-|, S$

$$
\xymatrix{
|-|:SS_{\ast}^{\mathcal{M}}\ar@<2pt>[r]&Top_{\ast}^{M}:S\ar@<2pt>[l]
}
.$$

 Define a functor \[Q:SS_{\ast}^{\mathcal{M}}\to SS_{\ast}^{\mathcal{M}}\] by $QX=S\bar{Q}|X|$ and a natural transformation 
$q:1\to Q$ by the composition \[X\stackrel{\eta_X}{\to} S|X|\stackrel{S\bar{q}_{|X|}}{\to}S\bar{Q}|X| \]
where $\eta$ is the unit of the adjunction. The proof of the following lemma is left as an exercise. 
It proceeds by reduction to 2.1. using well known
 properties of the pair  $|-|, S$.

\noindent{\bf Lemma 2.2.}\begin{it} The pair  $Q,q$ satisfies (A1) and (A2). \end{it}

Denote the subcategories of $Top_{\ast}^{M}$  whose objects have underlying spaces which are cofibrant by 
$Top^{M}_{\ast c}$
and write $SS^{\mathcal{M}}_{\ast f}$ for the subcategory of  $SS_{\ast}^{M}$ whose objects have fibrant
 underlying simplicial sets.
Let $Ab$ be the category of  abelian groups. We may consider $A\in Ab$ as a topological $M$-algebra in the obvious way.
This defines an inclusion functor $i:Ab\to Top^{M}_{\ast c}$.  
The assignment $M\to SiM$ defines a functor $S_{ab}:Ab\to SS^{\mathcal{M}}_{\ast f}$ from
abelian groups to $SS_{\ast}^{\mathcal{M}}$.

\noindent{\bf Lemma 2.3.  }\begin{it} The functor $S_{ab}$ is right adjoint to $\pi_0 Q$.

\end{it} 

\noindent{\bf Proof:} First note that by the adjunction  between simplicial and topological algebras there is a bijection:
\[Hom_{SS_{\ast}^{\mathcal{M}}}(X,S_{ab}A )\to Hom_{Top_{\ast}^M}(|X|,A).\]  

For any topological $M$-algebra $Y$ the map $Y\to  \pi_0 Y$ where $\pi_0 Y$ carries the quotient topology is
a morphism of $M$-algebras. 
 In case the space underlying $Y$ is of the homotopy type of a CW-complex the topology on 
$\pi_0 Y$ is the discrete one. This holds since the path components of a CW-complex are open and closed.
In particular this applies to an algebra which is cofibrant as a space because a generalized CW-complex 
is homotopy equivalent
 to a genuine CW-complex by cellular approximation. Now $|X|$ and $\bar{Q}|X|$ are cofibrant spaces.  
 It follows  that any morphism
$|X|\to A$ factorizes uniquely  over $|X| \to \pi_0|X|\to \pi_0 \bar{Q}|X|$ and this is the claim.
\hfill$\Box$

We need the following fact whose statement and proof are parallel to the ones of 
 \cite[5.4.]{BF}. 

\noindent{\bf Lemma 2.4.  }\begin{it} Every morphism $f:X\to Y$ in  $SS_{\ast}^ {\mathcal{M}}$   can be factored as 
\[X\stackrel{u}{\to}Z\stackrel{v}{\to}Y\] where $\pi_0 Qu:\pi_0 QX\to \pi_0 QZ$ is onto
 and $v$ is a $Q$-fibration. 
\end{it} 

\noindent{\bf Proof:} We define inductively a descending filtration of $\mathcal{M}$-algebras 
\[Y=C^0 \supset C^1 \supset \ldots C^{\alpha}\supset \ldots \]
indexed by the ordinals such that $f(X)\subset C^{\alpha}$ and $C^{\alpha}\subset Y$ is a $Q$-fibration as follows. 
Suppose $C^{\alpha}$ is found  define $C^{\alpha + 1}$ as the pullback 
$$
\xymatrix{
C^{\alpha +1} \ar[r] \ar[d] & 
S_{ab}G^{\alpha} \ar[d]
\\
C^{\alpha} \ar[r] & S_{ab}\pi_{0}\bar{Q}|C^{\alpha}|
}
$$ 
Where $G^{\alpha}$ is the image of $\pi_0 Qf:\pi_{0}QX\to \pi_{0}QC^{\alpha}$ and the map at the bottom is
 the composite \[C^{\alpha}\to S|C^{\alpha}| \to S\bar{Q}|C^{\alpha}|\to S_{ab}\pi_{0}\bar{Q}|C^{\alpha}|.\]
We claim that $S_{ab}G^{\alpha}\subset S_{ab}\pi_{0}QC^{\alpha}$ is a $Q$-fibration and hence $C^{\alpha + 1}\subset C^{\alpha}$
 is one as well. To see this let 
$$
\xymatrix{
A \ar[r]^{} \ar[d]_{} &  
S_{ab}G^{\alpha}{}\ar[d]^{}
\\
B \ar[r]_{}& S_{ab}\pi_{0}QC^{\alpha}
}
$$ 
be a commuting diagram with the vertical map on the left a $Q$-trivial cofibration.
   It is enough to  show existence of a filler for the adjoint diagram

$$
\xymatrix{
\pi_{0}Q A \ar[r]^{} \ar[d]_{} &  
G^{\alpha}{}\ar[d]^{}
\\
\pi_{0}Q B \ar[r]_{}& \pi_{0}QC^{\alpha}
}
$$

 The morphism on the left is an isomorphism between discrete abelian groups.
Hence the searched for filler exists.
Note that $f(X)\subset C^{\alpha + 1}$ and that $C^{\alpha + 1}\to Y$ is a fibration since fibrations
 are closed under composition.
 For a limit ordinal $\lambda$ such that $f(X)\subset C^{\alpha}$ for all $\alpha < \lambda$ define
 $C^{\lambda}=\lim_{\alpha < \lambda}C^{\alpha}$.
For sufficiently large $\alpha$ one has $C^{\alpha}=C^{\alpha + 1}$   and then  
$\pi_0 Qf:\pi_{0}QX\to \pi_{0}QC^{\alpha}$
is onto. Now put $Z=C^{\alpha}$.\hfill$\Box$

The next proposition establishes a weak form of axiom (A3) for the pair $Q,
q$.

\noindent{\bf Proposition 2.5. }\begin{it} Let 
$$
\xymatrix{
V \ar[r]^{k} \ar[d]_g & 
X \ar[d]^f
\\
W \ar[r]_h & Y
}
$$ 
be a pullback square of  $\mathcal{M}$-algebras with $f$ a fibration
 of fibrant objects such that
   the maps $q:X\to QX, q:Y\to QY$ and $Qh:QW\to QY$
 are  weak equivalences and  with $\pi_0 Qf:\pi_0 QX\to \pi_0 QY$ onto. Then 
$Qk:QV\to QX$ is a weak equivalence. 
 \end{it} 

\noindent{\bf Proof:}  Since $SG|Z|$ and $QZ$ are naturally weakly equivalent simplicial sets
 we may replace $\bar{Q}$ by $G$  in the argument.
 Consider the square of bisimplicial sets which results by application of the 
 functor \[Z\to SB_{\ast}G|Z|\] to the square above
and note that $SB_{\ast}=B_{\ast}S$ since $S$ commutes with products. Here $B_{\ast}X$ is the 
simplicial bar construction for topological or simplicial monoids $X$ which has the powers  $X^n$ in  simplicial degree $n$. 
  Now we wish to apply  \cite[B.4.]{BF}. This theorem  gives conditions on a square of bisimplicial sets which imply that the realization of the 
square is a homotopy pullback square in simplicial sets.\\
 We have to check that these conditions  hold in our situation.

First, we need to verify that the squares defined by $SG|X|^m, SG|Y|^m, 
\linebreak
SG|V|^m,$ $SG|W|^m$ are homotopy fibre squares. This is the case
since these spaces are naturally weakly equivalent to $X^m, Y^m, V^m, W^m$ and these form a fibre square
since $X,Y,V,W$ form one by assumption. 

Second, we need that the so called
$\pi_{\ast}$-Kan condition \cite[B.3.]{BF}  holds for $B_{\ast}SG|X|,B_{\ast}SG|Y|$ and that the map induced by $f$  
\[\alpha :\pi_{0}^v B_{\ast}SG|X|\to \pi_{0}^v B_{\ast}SG|Y|\] 
is a fibration. Here $\pi_{\ast}^{v}(Z)$ denotes the vertical homotopy of a bisimplicial set $Z$.
For simplicial spaces $Z$ which are simple in each degree the $\pi_{\ast}$-Kan condition is equivalent to the following condition \cite[B.3.1]{BF}: The obvious map
\[\beta :\pi_{t}^{v}Z_{free}\to \pi_{o}^{v}Z_{free}\]
is a fibration for each $t\geq 1$. Where for a simplicial set $U$, the symbol $\pi_{t}U_{free}$ stands for the set of  unpointed homotopy classes from $S^t$ to $|U|$.
By  assumption, $|X|, |Y|$ are group-like H-spaces. It follows that $B_{\ast}SX,B_{\ast}SY$
are indeed simple in each degree. 
Moreover, it follows that all the maps in question are surjective homomorphism of simplicial groups and hence
 are fibrations.\hfill$\Box$

A cofibrant fibrant approximation for a map $f:X\to Y$ in a model category is a commuting diagram

$$
\xymatrix{
X \ar[r]^{u} \ar[d]_{f} & 
 \widehat{X}\ar[d]^{ \hat{f}}
\\
Y \ar[r]_{v}& \widehat{Y}
}
$$ 
with trivial cofibrations $u,v$ and $\widehat{X}, \widehat{Y}$ fibrant. It is true that cofibrant fibrant approximation always exist and one can
choose  $ \widehat{f}$ as a fibration  \cite[8.1.23.]{HI}

Before we proceed we have to recall one more result of Bousfield. The proof of \cite[9.3.]{BO} gives us:

\noindent{\bf Proposition 2.6.}\begin{it} Let $\mathcal{D}$ be a proper model category and $Q,q$ a coaugmented functor
 for which (A1)-(A2) holds. Moreover, let $f:X\to Y$ be a map in $\mathcal{D}$   such that (A3) holds for one (and hence any)
 fibration  $\widehat{Qf}$ in a cofibrant fibrant approximation of $Qf$. Then $f$ can be factored into $ji$ with trivial
 $Q$-cofibration $i$ and $Q$-fibration $j$.  \end{it}

We are now ready to prove the main result of the paper:

\noindent{\bf Theorem 2.7.}\begin{it} Let $\mathcal{M}$ and $Q$ be as in 1.6. The category $SS_{\ast}^{\mathcal{M}}$ with the $Q$-structure is a left proper simplicial model category. Moreover, a morphism $f:X\to Y$ in $SS_{\ast}^{\mathcal{M}}$  is a $Q$-fibration if it is a fibration and 
$$
\xymatrix{
X \ar[r]^{q_{X}} \ar[d]_f & 
QX \ar[d]^{Qf}
\\
Y \ar[r]_{q_{Y}} & QY
}
$$ is a homotopy fibre square. In case $\pi_0 Qf$ is onto this condition is also necessary.
\end{it} 

\noindent{\bf Proof:} Since (A1) holds by 2.2. the axioms of a model category are satisfied by \cite[A.8]{BF} 
 except maybe  the trivial  cofibration fibration part of the factorization axiom.\\
Let $f:X\to Y$ be a map of $\mathcal{M}$-algebras. Factor $f=vu$ as in 2.4. in $Q$-fibration $v$ and with
$\pi_0 Qu$  onto. It is enough to factor $u$ into trivial $Q$-cofibration and $Q$-fibration. By 2.5.
and 2.6 this can be done.\\
To verify the left properness let 
$$
\xymatrix{
V \ar[r]^{k} \ar[d]_i & 
X \ar[d]^j
\\
W \ar[r]_h & Y
}
$$ be a pushout diagram with $Q$-equivalence $k$ and $Q$-cofibration $i$.
Factor $k$ in\-to $k=fg$ with $Q$-cofibration $g$ and trivial $Q$-fibration $f$.
Then $g$ is a trivial $Q$-cofibration and  $f$ is a trivial fibration in the underlying
 model category. The former by definition the later by \cite[A.8.(ii)]{BF}.
That the pushout of $g$ along $i$ is a $Q$-equivalence follows directly from the axioms of a model category.
 The pushout of $f$ along the induced cofibration $\bar{i}$ is a weak equivalence and hence a $Q$-equivalence
by the left properness of the underlying model category. The model structure on $SS_{\ast}^{\mathcal{M}}$ is simplicial by \cite{SV}. Hence, the $Q$-strucure is simplicial as well by \cite[9.7.]{BO} whose proof does not use (A3). The sufficiency of the stated condition follows from \cite[A.9.]{BF}. For the last statement, note that (A3) holds for $f$ by 2.5. Now the proof proceeds as in \cite[A.10.]{BF}.
\hfill$\Box$

\noindent{\bf Corollary 2.8.}\begin{it} An object $X\in SS_{\ast}^{\mathcal{M}}$ is fibrant in the $Q$-structure if and only if the underlying simplicial set is fibrant and $q_X$ is a weak equivalence.\end{it}

\noindent{\bf Remark 2.9.} The $Q$-structure on $SS_{\ast}^{\mathcal{M}}$ is not right proper.
This follows from the fact that (A3) does not hold in general. An example which shows this can be found in \cite{BF}
on page 109. 

\noindent{\bf Remark 2.10.} One can show that the $\bar{Q}$-structure induced  on $M$-algebras satisfies axioms which are slightly weaker than those of a cofibration category but still strong enough to induce a well defined homotopy category. Most of the axioms hold only if the source (and sometimes the target) of the morphisms are cofibrant spaces. One has to use the modifications for some of the arguments in \cite{BF} which were already hinted on by Bousfield in \cite[9.5.]{BO}.

\begin{it}
\begin{center}\noindent M. Stelzer\\
Universit\"at Osnabr\"uck\\
Fachbereich Mathematik/Informatik\\ 
 49069 Osnabr\"uck\\
Germany\\
\noindent
Email address: mstelzer@uos.de\\
\end{center}\end{it}

\begin{thebibliography}{99999}
\bibitem{BCV} B.Badzioch, K. Chung, A. Voronov, The cannonical delooping machine, Journal of Pure and Applied Algebra 208(2) (2007) 531-540.
\bibitem{B}C. Balteanu, Z. Fiedorowicz, R. Schwaenzl, R.M.  Vogt, Iterated monoidal categories, Advances in Mathematics 176 (2003) 277-349.
\bibitem{Ba}H. Baues, Algebraic Homotopy, Cambridge studies in advanced mathematics 15 (1989).
\bibitem{BM1}C. Berger, I. Moerdijk, Axiomatic homotopy theory for operads, Comment. Math. Helv. 78 (2003) 805-831.
\bibitem{BM2}C. Berger, I. Moerdijk, The Boarman-Vogt Resolution Of Operads In Monoidal Model Categories, Topology 45 (2006) 807-849.
\bibitem{BV}J.M. Boardman, R.M.  Vogt, Homotopy Invariant Algebraic Structures on Topological Spaces, Lecture Notes in Math. 347, Springer-Verlag, (1973).
\bibitem{BO} A.K. Bousfield, On the Telescopic Homotopy Theory of Spaces, Trans. Amer. Math. Soc. 353, Number 6, (2001) 2391-2426.
\bibitem{BO2} A.K. Bousfield, The simplicial homotopy theory of iterated loop spaces, Manuscript,(1992). 
\bibitem{BF} A.K. Bousfield, E.M. Friedlander, Homotopy theory of $\Gamma$-spaces, Spectra, And Bisimplicial Sets, Geometric applications of homotopy theory 2, Lecture Notes in Math. 658, Springer Verlag (1977).
\bibitem{BFV}M. Brun, Z. Fiedorowicz, R. Vogt, On the multiplicative structure of Hochschild homology, Algebraic and Geometric Topology 7, 1633-1650 (2007).
\bibitem{C}F. Cohen, The homology of $\mathcal{C}_{n+1}$ spaces $n\geq 0$, Lecture Notes in Math. 533, Springer-Verlag, (1976), 207-351. 
\bibitem{Cr}S.E. Crans, Quillen closed model structures for sheaves, Journal of Pure and  Applied Algebra 101 (1995) 35-57.
\bibitem{FV}Z. Fiedorowicz, R.M. Vogt, Simplicial $n$-Fold Monoidal Categories Model All Loop Spaces, Cahiers de Topologie et Geometrie Differentielle Categorique, Volume XLIV-2 (2003) 105-148.
\bibitem{FSV}Z. Fiedorowicz, M. Stelzer, R.M. Vogt, Homotopy Colimits of Al\-ge\-bras over Cat-Operads and Iterated Loop Spaces, arXive:math/1109.0265.
\bibitem{GJ}P. Goerss, R. Jardine, Simplicial homotopy theory, Birkhaeuser Verlag, Basel (1999).
\bibitem{HI}P.S. Hirschhorn, Model Categories and Their Localizations, Mathematical Surveys and Monographs, Vol. 99. Amer. Mathc., Providence,RI, (2002).
\bibitem{H}M. Hovey, Model Categories, Mathematical Surveys and Monographs, Vol. 63. Amer. Mathc., Providence,RI, (1999).
 Series 64 (1982).
\bibitem{MSS}M. Markl, S. Shnider, J. Stasheff, Operads in  Algebra, Topology and Physics, Mathematical Surveys and Monographs, Vol. 96. Amer. Math., Providence,RI, (2002).
\bibitem{M}J.P. May, The Geometry of Iterated  Loop Spaces, Lecture Notes in Math. 271, Springer-Verlag, (1972).
\bibitem{M2}J.P. May, $E_{\infty}$ spaces, group completions, and permutative categories, London Math. Soc. Lecture Notes No. 11 (1974), 61-93.
\bibitem{SMS}D. McDuff, G.Segal, Homology Fibrations and the Group-Com\-pletion Theorem, Inventiones Math. 31 (1976) 279-284.
\bibitem{MI}R.J. Milgram, The Bar Construction and Abelian H-Spaces, Illinois Journal of Mathematics 11 (1967) 242-250.
\bibitem{Q}D.G. Quillen, Homotopical Algebra, Lecture Notes in Math. 43, Springer Verlag, (1967).
\bibitem{QFM}D.G. Quillen, E.M. Friedlander, B. Mazur, Filtrations on the homology of algebraic varieties, with an appendix by D. Quillen, Mem. Amer. Math. Soc. 110(529), (1994) x+110.
\bibitem{S}G. Segal, Categories and Cohomology Theories, Topology 13 (1974) 293-312.
\bibitem{S2}G.Segal, Configuration spaces and iterated loop spaces, Inventiones Math. 21 (1973), 213-221.
\bibitem{SV}R. Schw\"anzl, R.M. Vogt, The categories of $A_{\infty}-$ and $E_{\infty}-$monoids and ring spaces as closed simplicial and topological model categories, Arch. Math. 56 (1991) 405-411.
\bibitem{SP} M. Spitzweck, Operads, algebras and modules in general model categories, PhD thesis, Bonn (2001) arXive:math/0101102.
\bibitem{STa} A.E. Stanculescu, Note on a theorem of Bousfield and Friedlander, Topology and its Applications 155, (2008) 1434-1438.
\bibitem{ST} A. Str\o m, The homotopy category is a homotopy category, Arch. Math. 22 (1972) 435-441.
\bibitem{T}R.W. Thomason, Uniqueness of Delooping Machines, Duke Mathematical Journal 46, No.2 (1976) 217-252.
\bibitem{V1}R.M. Vogt, Convenient categories of topological spaces for homotopy theory, Arch.  Math. 22 (1971) 545-555.
\bibitem{V3}R.M. Vogt, Cofibrant operads and universal $E_{\infty}$-operads, Topology and its Applications 133 (2003) 69-87.

\end{thebibliography}
 \end{document}